\documentclass[10pt, leqno]{amsart}
\usepackage{amssymb,amsfonts,amscd,amsbsy}
\usepackage[mathscr]{eucal}

\theoremstyle{plain}
\newtheorem{thm}{Theorem}[section]
\newtheorem{lem}[thm]{Lemma}
\newtheorem{pro}[thm]{Proposition}
\newtheorem{cor}[thm]{Corollary}
\theoremstyle{definition}

\newtheorem{rem}[thm]{Remark}
\newtheorem{exa}[thm]{Example}

\begin{document}
\title{Annihilating polynomials for quadratic forms and Stirling numbers of the second kind}
\author{Stefan De Wannemacker}
\maketitle

\begin{abstract}
We present a set of generators of the full annihilator ideal for
the Witt ring of an arbitrary field of characteristic unequal to
two satisfying a non-vanishing condition on the powers of the
fundamental ideal in the torsion part of the Witt ring.  This
settles a conjecture of Ongenae and Van Geel.  This result could
only be proved by first obtaining a new lower bound on the 2-adic
valuation of Stirling numbers of the second kind.
\end{abstract}
%%%%%%%%%%%%%%%%%%%%%%%%%%%%%%%%%%%%%%%%%%%%%%%%%%%%%%%%%%%%%%%%%%%%%%%%%%%%%%%%%%%%%%%%%%%%%%%%%%%%%%%%%%%%%%%%%%
%                                                                                                                %
%  SECTION : INTRODUCTION                                                                                        %
%                                                                                                                %
%%%%%%%%%%%%%%%%%%%%%%%%%%%%%%%%%%%%%%%%%%%%%%%%%%%%%%%%%%%%%%%%%%%%%%%%%%%%%%%%%%%%%%%%%%%%%%%%%%%%%%%%%%%%%%%%%%
\section{Introduction}

In 1937, Witt already observed that the  Witt ring was integral in
the sense that each element was annihilated by a monic integer
polynomial.  Fifty years later, in 1987, Lewis was the first to
give explicit examples of such polynomials \cite{LEW2}. He showed
that the monic polynomial, $p_n(X)$, defined as
\[p_n(X)=(X-n)(X-(n-2))\ldots (X+(n-2))(X+n)\]
annihilates every non-singular quadratic form of dimension $n$
over every field $F$ of characteristic unequal to two.\\ Since
then, many other polynomials in $\mathbb{Z}[X]$ were found
annihilating all or a family of classes of nonsingular quadratic
forms in the Witt ring and we refer the reader to \cite{LEW} for a
nice survey of the main results on this topic.\\

Let $F$ be a field of characteristic not $2$. The object we want
to consider here is the torsion annihilator ideal
\[ A_t(F) = \{f(X) \in \mathbb{Z}[X]\ |\ f(\phi)=0,\ \forall\ \phi \in I_t(F)\}\]
where $I_t(F)=W_t(F)\,\cap\, I(F)$, $\ W_t(F)$ the torsion part of
the Witt ring and $I(F)$ the ideal of all even-dimensional forms
in the Witt ring. Since $A_t(F)$ is an ideal in the noetherian
ring $\mathbb{Z}[X]$, it is finitely generated. The main problem
is to find a set of generators for this ideal.
\ \\
We will prove the following result.\\
\ \\
\emph{For fields $F$ satisfying the conditions that $2^rW_t(F)=0$
and $2^{r-1}(I_t(F))^{2k-1}\neq 0$ with $k$ uniquely determined by
$r$, the torsion annihilator ideal $A_t(F)$ is the ideal generated
by the monomials\\
\[{\large \lbrace 2^r \rbrace \cup \lbrace 2^{r-\nu_2((2i)!)}X^{2i}\rbrace_{1 \leq i \leq
k-1}\ \cup\ \lbrace X^{2k}\rbrace\qquad}\text{for nonreal fields
and}\]
\[\lbrace 2^rX \rbrace \cup \lbrace 2^{r-\nu_2((2i)!)}X^{2i}\rbrace_{1 \leq i \leq
k-1}\ \cup\ \lbrace X^{2k}\rbrace\qquad\text{for real fields, }\]
where $\nu_2$ denotes the $2$-adic valuation.
}\\
\ \\
In the case of a nonreal field $F$ this theorem was conjectured (see
Corollary \ref{coro_conj} ) by Ongenae and Van Geel \cite{OVG}. They
gave a proof for fields with level $s(F) \leq 16$ and using the same
technique, one can check that the theorem holds for all nonreal
fields $F$ with level $s(F)\leq 64$, but a general method was
lacking.
\ \\
The general method, used to prove the theorem, consists in
evaluating a polynomial $f(X) \in \mathbb{Z}[X]$ in the
even-dimensional forms $\bot_{i=1}^{n}\langle\!\langle
a_i\rangle\!\rangle$ with $1 \leq n \leq deg(f).$  This evaluation
can be rewritten as a linear combination of sums of Pfister forms
and the coefficients that appear turn out to be related to the
Stirling numbers of the second kind.\
\ \\
The result about the (torsion) annihilator ideal could only be
proved by first obtaining a new lower bound for the $2$-adic
valuation of all Stirling numbers $S(n,k)$ of the second kind,
namely
\[\nu_2(S(n,k))\geq d(k)-d(n),\quad\text{for}\quad0< k\leq n\]
where $d(k)$ is the sum of the binary digits in the binary representation of $k$.\\
%%%%%%%%%%%%%%%%%%%%%%%%%%%%%%%%%%%%%%%%%%%%%%%%%%%%%%%%%%%%%%%%%%%%%%%%%%%%%%%%%%%%%%%%%%%%%%%%%%%%%%%%%%%%%%%%%%
%                                                                                                                %
%  SUBSECTION : COMBINATORICS                                                                                    %
%                                                                                                                %
%%%%%%%%%%%%%%%%%%%%%%%%%%%%%%%%%%%%%%%%%%%%%%%%%%%%%%%%%%%%%%%%%%%%%%%%%%%%%%%%%%%%%%%%%%%%%%%%%%%%%%%%%%%%%%%%%%
\section{Stirling numbers of the second kind}
\subsection{Preliminaries}

Let $n \in \mathbb{N}\,.$ The Stirling numbers $S(n,k)\ (k \in
\mathbb{N})$ of the second kind are given by
\[x^n=\sum_{k=0}^{\infty}S(n,k)(x)_k,\]
where $(x)_k=x(x-1)(x-2)\ldots(x-k+1)$ for $k \in
\mathbb{N}\setminus \{0\}$ and $(x)_0=1.$  Actually $S(n,k)$ is
the number of ways in which it is possible to partition a set with
$n$ elements in $k$\ classes.\\
The Stirling numbers of the second kind can be computed in several
ways.
\begin{pro}\label{prop}
\begin{align*}
S(n,k) &= \frac{1}{k!} \sum_{i=0}^{k}(-1)^i{k\choose k-i}(k-i)^n,\\
S(n,k)&=S(n-1,k-1)+kS(n-1,k)\\
&\text{ with }\ S(n,0)=S(0,k)=0\ \text{ and }\ S(0,0)=1\\
S(n,k)&=\frac{1}{k!} \sum_{n_1,n_2,...,n_k}{n\choose n_1,n_2,...n_k},\\
&\text{where } n_1,n_2,...,n_k \text{ are non-zero and their sum equals } n.\\
\end{align*}
\end{pro}
\begin{proof}
\ \\
See \cite{KN} and \cite{R}.\\
\end{proof}

\subsection{$2$-adic valuation of Stirling numbers of the second kind}
The $2$-adic valuation of Stirling numbers of the second kind and
the $2$-adic valuations of other combinatorial numbers have been
widely studied, but many problems in this area are still unsolved.
We will give a new lower bound for the $2$-adic valuation of all
Stirling numbers of the second kind.\\
\ \\
Denote by $d(n)$ the sum of the digits in the binary
representation of $n$ and define the $2$-adic valuation function
$\nu_2(n)$ for all non-zero integers $n$ by $\nu_2(n)=p$, where
$2^p|n$ and $2^{p+1}\nmid n$.\\
Recall the following properties (see \cite{KN}).
\begin{align*}
\nu_2(n!)&=n-d(n)\qquad\text{ (Legendre)}\\
\nu_2\left(\binom{n}{k}\right)&=d(k)+d(n-k)-d(n)\qquad\text{ (Kummer)}\\
\end{align*}
for all $k,n \in \mathbb{N}$ with $0\leq k\leq n$.\\
A new lower bound on the $2$-adic valuation of Stirling numbers of
the second kind can be obtained as follows.\\
\begin{cor}\label{new_lower_bound} Let $n,k \in \mathbb{N}$ and $0 < k \leq
n$. Then
\[\nu_2\left(S(n,k)\right) \geq d(k) - d(n).\]
\end{cor}
\begin{proof}
$\text{ }$\\
\ \\
Using the fact that $\nu_2(a+b)\geq {\rm min}\{\nu_2(a),\nu_2(b)\}$,
the result follows from  \ref{prop} and
\begin{align*}
\nu_2\left(\frac{1}{k!}\binom{n}{n_1,\ldots,n_k}\right)&=\nu_2(n!)-\nu_2(k!)-\sum_{i=1}^{k}\nu_2(n_i!)\\
&=n-d(n)-(k-d(k))-(\sum_{i=1}^k(n_i-d(n_i)))\\
& \text{\qquad\qquad\qquad\qquad\qquad\qquad(by the Legendre identity)}\\
&\geq d(k)-d(n)\\
& \text{\qquad\qquad\qquad\qquad\qquad\qquad(since $d(n_i)\geq 1$).}\\
\end{align*}
\end{proof}
\subsection{Relationship between Stirling numbers of the second kind and quadratic forms}
\ \\
Consider the $\mathbb{Z}-$algebra $B =
\mathbb{Z}[X_1,X_2,\ldots,X_n]/K$ where
$K=(X_1^2-2X_1,X_2^2-2X_2,\ldots,X_n^2-2X_n)$. This algebra has a
$\mathbb{Z}-$basis given by the elements $Y_J :=\prod_{j\in J}X_j
\mod K$ where $\emptyset\neq J \subseteq \{1,2,\ldots,n\}$ and
$Y_{\emptyset}:=1$.\\
Stirling numbers of the second kind turn up in a natural way by
making calculations in this $\mathbb{Z}-$algebra.  The calculations
frequently use the following relation
\[Y_{\{i\}}^p=2^{p-1}Y_{\{i\}}.\]
 We have the following
\begin{pro}\label{zalgebraprop} Let $f(X)=c_dX^d+\ldots+c_1X + c_0\in
\mathbb{Z}[X].$  Then
\begin{align*}
f\left(\sum_{i=1}^{n}Y_{\{i\}}\right)&=\ \left(\sum_{q=0}^d2^{q-0}0!\ S(q,0)c_q\right)Y_{\emptyset}\\
&\\
&+\ \left(\sum_{q=1}^d2^{q-1}1!\ S(q,1)c_q\right)\left(\sum_{i=1}^nY_{\{i\}}\right)\\
&\\
&+\ \left(\sum_{q=2}^d2^{q-2}2!\ S(q,2)c_q\right)\left(\sum_{i<j}^{n}Y_{\{i,j\}}\right)\\
&+\ \ldots\\
&+\ \left(\sum_{q=n}^d2^{q-n}n!\
S(q,n)c_q\right)Y_{\{1,2,\ldots,n\}}.
\end{align*}
\end{pro}
\begin{proof}
The evaluation of the polynomial $f(X)=c_dX^d+\ldots+c_1X+c_0$ in
sums of the basis elements $\sum_{i=1}^{n}Y_{\{i\}}$ can be written
in the following unique way.
\begin{align*}
f(\sum_{i=1}^{n}Y_{\{i\}}) &=\ A_0(c_0,\ldots ,c_d)Y_{\emptyset}\\
&+\ A_1(c_0,\ldots ,c_d)(\sum_{i=1}^nY_{\{i\}})\\
&+\ A_2(c_0,\ldots ,c_d)(\sum_{i<j}^{n}Y_{\{i,j\}})\\
&+\ \ldots\\
&+\ A_{n}(c_0,\ldots ,c_d)Y_{\{1,2,\ldots,n\}},
\end{align*}
where
\[A_p(c_0,\ldots,c_d)=\sum_{q=0}^d\gamma_{p,q}c_q\]
with $\gamma_{p,q}\in\mathbb{N}$.
For $p=0$, $A_0(c_0,\ldots,c_d)=c_0=\sum_{q=0}^d2^{q-0}0!\ S(q,0)c_q.$\\
For $p\geq 1$, $\gamma_{p,q}$ is the coefficient of
$Y_{\{1,2,\ldots,p\}}$ in
$(\sum_{i=1}^{p}Y_{\{i\}})^{q}$. For $1 \leq p\leq q$ we can write\\
\begin{align*}
\gamma_{p,q}&=\sum_{q_1,q_2,...,q_p \geq 1}^{q_1+q_2+...+q_p=q}{q\choose q_1,q_2,...,q_p}2^{q_1-1}2^{q_2-1}...\ 2^{q_p-1}\\
&=2^{q-p}\sum_{q_1,q_2,...,q_p \geq 1}^{q_1+q_2+...+q_p=q}{q\choose q_1,q_2,...,q_p}\\
&=2^{q-p}p!\,S(q,p).
\end{align*}\\
If $p>q$ then clearly no basis-element $Y_{\{1,2,\ldots,p\}}$ can
occur in $(\sum_{i=1}^{p}Y_{\{i\}})^{q}$. So,
\[\gamma_{p,q}=
\begin{cases} 2^{q-p}p!\,S(q,p) & \qquad\text{if} \quad p \leq q,\\
0 &\qquad\text{otherwise.} \end{cases}
\]
Applying this to the coefficients $A_p(c_0,\ldots,c_d)$ we obtain
\begin{align*}
f\left(\sum_{i=1}^{n}Y_{\{i\}}\right)&=\ \left(\sum_{q=0}^d2^{q-0}0!\ S(q,0)c_q\right)Y_{\emptyset}\\
&\\
&+\ \left(\sum_{q=1}^d2^{q-1}1!\ S(q,1)c_q\right)\left(\sum_{i=1}^nY_{\{i\}}\right)\\
&+\ \left(\sum_{q=2}^d2^{q-2}2!\ S(q,2)c_q\right)\left(\sum_{i<j}^{n}Y_{\{i,j\}}\right)\\
&+\ \ldots\\
&+\ \left(\sum_{q=n}^d2^{q-n}n!\ S(q,n)c_q\right)Y_{\{1,2,\ldots,n\}}.\\
\end{align*}
\end{proof}
We can evaluate \\
\[f(X)=c_dX^d+\ldots+c_1X \in \mathbb{Z}[X]\]\\
in classes of quadratic forms $\phi \in W(F)$ by defining\\
\[f(\phi)=c_d\phi^d\ \bot\ \ldots\ \bot\ c_1\phi\ \,\in W(F)\]\\
where $c_i\phi = sign(c_i)\underbrace{(\phi\ \bot\ \ldots\ \bot\
\phi)}_{|c_i| \ \text{times}}$ and $\phi^i = \underbrace{\phi\
\otimes\ \ldots\
\otimes\ \phi}_{i \ \text{times}}$.\\
\ \\
For arbitrary $k>0$ and $\ a_1,a_2,\ldots,a_k \in
F^{\ast}=F\setminus\{0\}$, let $\langle\!\langle
a_1,a_2,\ldots,a_k\rangle\!\rangle$ denote the $2^k$-dimensional
$k$-fold Pfister form
\[\langle\!\langle
a_1,a_2,\ldots,a_k\rangle\!\rangle :=\langle 1,a_1\rangle \otimes
\langle 1,a_2\rangle \otimes \ldots \otimes \langle 1,a_k\rangle.\]\\
Observe that the following relation holds.
\[\langle\!\langle a \rangle\!\rangle^p = 2^{p-1} \langle\!\langle a
\rangle\!\rangle.\] \ \\
\ \\
\begin{cor}\label{eval_proposition} Let $f(X)=c_dX^d+\ldots+c_1X + c_0 \in \mathbb{Z}[X].$\\
Then
\begin{align*}
f(\bot_{i=1}^{n}\langle\!\langle a_i\rangle\!\rangle)&=\ \left(\sum_{q=0}^d2^{q-0}0!\ S(q,0)c_q\right)\langle 1 \rangle\\
&\bot\ \left(\sum_{q=1}^d2^{q-1}1!\ S(q,1)c_q\right)(\bot_{i=1}^n\langle\!\langle a_i\rangle\!\rangle)\\
&\bot\ \left(\sum_{q=2}^d2^{q-2}2!\ S(q,2)c_q\right)(\bot_{i<j}^{n}\langle\!\langle a_i,a_j \rangle\!\rangle)\\
&\bot\ \ldots\\
&\bot\ \left(\sum_{q=n}^d2^{q-n}n!\ S(q,n)c_q\right)\langle\!\langle
a_1,\ldots,a_n\rangle\!\rangle.
\end{align*}
\end{cor}
\begin{proof} Since we have a $\mathbb{Z}-$algebra homomorphism $B \rightarrow W(F)$ with $Y_{\{i\}}\mapsto\langle\!\langle a_i
\rangle\!\rangle$ and $Y_{\{J\}}\mapsto\otimes_{j\in
J}\langle\!\langle a_i \rangle\!\rangle$ for all $J\subseteq
\{1,2,\ldots,n\}$ this result follows from the Proposition
\ref{zalgebraprop} by substituting $\sum_{i=1}^{n}Y_{\{i\}}$ by
$\bot_{i=1}^{n}\langle\!\langle a_i\rangle\!\rangle$.
\end{proof}
\begin{cor}\label{coro}
Let $f(X)=c_dX^d+\ldots+c_1X \in \mathbb{Z}[X]$ and $ \phi\ =\
\langle a_1,a_2,\ldots a_n\rangle$ a quadratic form of dimension
$n$. Then
\begin{align*}
f(\phi)&=\ \left(\sum_{q=0}^d\sum_{t=q}^d2^{q-0}0!\ S(q,0){t\choose q}(-n)^{t-q}c_t\right)\langle 1 \rangle\\
&\bot\ \left(\sum_{q=1}^d\sum_{t=q}^d2^{q-1}1!\ S(q,1){t\choose q}(-n)^{t-q}c_t\right)(\bot_{i=1}^n\langle\!\langle a_i\rangle\!\rangle)\\
&\bot\ \left(\sum_{q=2}^d\sum_{t=q}^d2^{q-2}2!\ S(q,2){t\choose q}(-n)^{t-q}c_t\right)(\bot_{i<j}^{n}\langle\!\langle a_i,a_j \rangle\!\rangle)\\
&\bot\ \ldots\\
&\bot\ \left(\sum_{q=n}^d\sum_{t=q}^d2^{q-n}n!\ S(q,n){t\choose
q}(-n)^{t-q}c_t\right)\langle\!\langle
a_1,\ldots,a_n\rangle\!\rangle
\end{align*}
\end{cor}
\begin{proof}
\ \\
Note that for all
\[ \phi\ =\ \langle a_1,a_2,\ldots a_n\rangle \]
 we have
\[ \phi\ =\ \bot_{i=1}^{n}\langle\!\langle
a_i\rangle\!\rangle - n\langle 1 \rangle. \] Using the Taylor series
of $f$, we have
\begin{align*}
f(X)&=\sum_{q=0}^\infty \frac{1}{q!} f^{(q)}(-n)(X+n)^q\\
&=\sum_{q=0}^d\sum_{t=q}^d{t\choose q}(-n)^{t-q}c_t(X+n)^q\\
&\\
&=:g(X+n).
\end{align*}
The result follows from the previous proposition applied to $g$ and
the $\mathbb{Z}-$algebra homomorphism $B \rightarrow W(F)$ described
in the proof of the previous corollary.
\end{proof}
\begin{lem}\label{equations}
Let $f(X)=c_dX^d+\ldots+c_1X \in \mathbb{Z}[X]$.  Let
$n\in\mathbb{N},\ n\leq d,\ a_1,a_2,\ldots ,a_d \in F^{\ast}$ . If
$\ f(\bot_{i=1}^{k}\langle\!\langle a_{\sigma
(i)}\rangle\!\rangle)=0,\ $ for all $\ 1\leq k \leq n,\
 \sigma \in S_d$ then
\[\left(\sum_{q=n}^{d}2^{q-n}n!\
S(q,n)c_q\right)\langle\!\langle a_{\tau(1)},\ldots a_{\tau(n)}
\rangle\!\rangle=0,\quad\text{for all}\quad \tau \in S_d.\]
\end{lem}
\begin{proof}
\ \\
By induction on $n$.\\
For $n=1$,
\begin{align*}
0&= f(\langle\!\langle a_i\rangle\!\rangle)\\
&= c_d\langle\!\langle
a_i\rangle\!\rangle^d+c_{d-1}\langle\!\langle
a_i\rangle\!\rangle^{d-1}+\ldots+c_1\langle\!\langle
a_i\rangle\!\rangle\\
&= (2^{d-1}c_d+2^{d-2}c_{d-1}+\ldots+c_1)\langle\!\langle
a_i\rangle\!\rangle\\
&= \left(\sum_{q=1}^d2^{q-1}1!\ S(q,1)c_q\right)\langle\!\langle
a_i\rangle\!\rangle
\end{align*}
since $S(q,1)=1$ for all $q\geq 1$.\\
Assume now that the lemma is true for all $i<n$.  We will
prove the lemma for $n$.\\
Let $\tau \in S_d$.
\begin{align*}
0=&\ f(\bot_{i=1}^{n}\langle\!\langle
a_{\tau(i)}\rangle\!\rangle)\\
&\\
=&\ \left(\sum_{q=1}^d2^{q-1}1!\ S(q,1)c_q\right)(\bot_{i=1}^n\langle\!\langle a_{\tau(i)}\rangle\!\rangle)\\
&\bot\ \left(\sum_{q=2}^d2^{q-2}2!\ S(q,2)c_q\right)(\bot_{\tau(i)<\tau(j)}^{n}\langle\!\langle a_{\tau(i)},a_{\tau(j)} \rangle\!\rangle)\\
&\bot\ \ldots\\
&\bot\ \left(\sum_{q=n}^d2^{q-n}n!\
S(q,n)c_q\right)\langle\!\langle
a_{\tau(1)},\ldots,a_{\tau(n)}\rangle\!\rangle\qquad (1)\\
\end{align*}
by Corollary \ref{eval_proposition}.\\
Let $1 \leq k\leq n-1.$  Note that for every subset $U
\subset\{1,\ldots,d\}$ of $k$ elements, there exists a permutation
$\sigma \in S_d$ such that
$U=\{\sigma(1),\ldots,\sigma(k)\}$.\\
So,
\begin{align*}
&\left(\sum_{q=k}^d2^{q-k}k!\ S(q,k)c_q\right)(\langle\!\langle a_{\tau(i_1)},\ldots,a_{\tau(i_k)} \rangle\!\rangle)\\
& = \left(\sum_{q=k}^d2^{q-k}k!\ S(q,k)c_q\right)(\langle\!\langle a_{\sigma(1)},\ldots,a_{\sigma(k)}\rangle\!\rangle)\\
& = 0\\
\end{align*}
by the induction hypothesis. So, $(1)$ becomes
\[
0=\ \left(\sum_{q=n}^d2^{q-n}n!\ S(q,n)c_q\right)\langle\!\langle
a_{\tau(1)},\ldots,a_{\tau(n)}\rangle\!\rangle\]
\end{proof}
%%%%%%%%%%%%%%%%%%%%%%%%%%%%%%%%%%%%%%%%%%%%%%%%%%%%%%%%%%%%%%%%%%%%%%%%%%%%%%%%%%%%%%%%%%%%%%%%%%%%%%%%%%%%%%%%%%
%                                                                                                                %
%  SECTION : Annihilator ideals                                                                                  %
%                                                                                                                %
%%%%%%%%%%%%%%%%%%%%%%%%%%%%%%%%%%%%%%%%%%%%%%%%%%%%%%%%%%%%%%%%%%%%%%%%%%%%%%%%%%%%%%%%%%%%%%%%%%%%%%%%%%%%%%%%%%
\section{Polynomials annihilating the Witt Ring}
\subsection{Preliminaries}
%\begin{lemma}\label{unique_k}
%For all $r \in \mathbb{N}$ there exists a unique $k \in
%\mathbb{N}$ such that
%\[\nu_2((2k-2)!)<r\leq \nu_2((2k)!).\]
%\end{lemma}
%\begin{proof}
%See \cite{OVG}.
%\end{proof}
\ \\
The \textit{fundamental ideal} $I(F)$ of $W(F)$ is the ideal
consisting of all even-dimensional forms of $W(F)$. Put
$I_t(F):=I(F)\cap W_t(F)$.\\
A field $F$ is called \textit{nonreal} if $-1$ is a sum of squares
in $F$, in which case we define the \textit{level} of $F$ to be
\[ s(F):= min\{n \in \mathbb{N} \;\vline\;
a_1^2+a_2^2+\ldots+a_n^2=-1,\ a_i \in F\}.\] Otherwise, $F$ is
called \textit{(formally) real} and the level is defined to be
$s(F)=\infty$.\\ It is a well-known fact that $s(F)$ is a power of
two, if $F$ is nonreal (\cite{SCH} or \cite{LAM}).  There is no odd
torsion in the Witt ring.  The \textit{height} of $F$ is defined to
be $h(F):=\infty$ if there is no $2$-power $2^d$ with
$2^dW_{t}(F)=0$.  Otherwise, $h(F)$ is the smallest such $2$-power
with $2^dW_{t}(F)=0$.\\
%For real fields $F$, $I_t(F)=W_t(F)$ since every torsion is even.
%For nonreal fields $F$, $I_t(F)=I(F)$ since $W_t(F) = W(F)$.
%We will need the following result.
%\begin{lemma}\label{existence_of_elements}
%Let $F$ be a field such that $(I_t(F))^n\neq0$.  Then there exist
%$a_1,\ldots,a_n \in F^{\ast}$ such that \begin{itemize} \item[(i)]
%$\langle\!\langle a_i \rangle\!\rangle\in I_t(F)$,\quad for all\
%$1 \leq i\leq n$,\item[(ii)] $\langle\!\langle
%a_{1},\ldots,a_{n}\rangle\!\rangle\neq0$.\end{itemize}
%\end{lemma}
%\begin{proof}
%For a nonreal field $F$, $I_t(F)=I(F)$. Since $I^n(F)$ is finitely
%generated by the $n$-fold Pfister forms, the non-vanishing
%condition on the power of the fundamental ideal implies the
%existence of elements $a_1,\ldots,a_{n} \in F^{\ast}$ such that
%$\langle\!\langle a_{1},\ldots,a_{n}\rangle\!\rangle\neq 0$.\\
%For a real field $F$, $I_t(F)=W_t(F)$. The torsion ideal $W_t(F)$
%is generated by the $1$-fold Pfister forms $\langle\!\langle a
%\rangle\!\rangle$ where $a$ is totally negative.  The condition on
%the power of this torsion ideal thus implies the existence of
%totally negative elements $a_{1},\ldots,a_{n}$ such that
%$\langle\!\langle a_{1},\ldots,a_{n}\rangle\!\rangle\neq 0$.
%\end{proof}
 \vskip 15pt \noindent We define the \textit{torsion annihilator
ideal} $A_{t}(F)$ in $\mathbb{Z}[X]$ by
\[ A_{t}(F)=\{ f(X)\in \mathbb{Z}[X]\;\vline\; f(\varphi)=0 \;\text{for all}\; \varphi \in I_{t}(F)\}.\]
\ \\
For $F$ nonreal, define the \textit{full annihilator ideal} $A(F)$
in $\mathbb{Z}[X]$ by
\[ A(F)=\{ f(X)\in \mathbb{Z}[X]\;\vline\; f(\varphi)=0 \;\text{for all}\; \varphi \in W(F)\},\]
the \textit{even annihilator ideal} by\\
 \[A_e(F)=\{f(X)\in \mathbb{Z}[X]\; \vline\; f(\varphi)=0 \;\text{for all even-dimensional}\; \varphi \in W(F)\}\]
and the \textit{odd annihilator ideal} by\\
\[ A_o(F)=\{ f(X)\in \mathbb{Z}[X]\; \vline\; f(\varphi)=0 \;\text{for all odd-dimensional}\;\varphi \in W(F)\}.\]
\ \\
In what follows, let $k=k(r)$ be the natural number uniquely
determined by $\nu_2((2k-2)!) < r \leq \nu_2((2k)!)$ (see
\cite{OVG}).
\ \\
\ \\ We define the ideals\\
\begin{eqnarray*}
J'_{e,r} &=& (2^{r}X)+(\{2^{r-\nu_2((2i)!)}X^{2i}\}_{1\leq i\leq k-1})\ +\  (X^{2k})\\
&=& (2^{r}X,2^{r-1}X^{2},2^{r-3}X^{4},...,2^{r-\nu_2((2i)!)}X^{2i},...,2^{r-\nu_2((2k-2)!)}X^{2k-2},X^{2k}),\\
J_{e,r} &=& (\{2^{r-\nu_2((2i)!)}X^{2i}\}_{0\leq i\leq k-1})\ +\  (X^{2k})\\
&=& (2^{r},2^{r-1}X^{2},2^{r-3}X^{4},...,2^{r-\nu_2((2i)!)}X^{2i},...,2^{r-\nu_2((2k-2)!)}X^{2k-2},X^{2k}),\\
J_{o,r}&=& (\{2^{r-\nu_2((2i)!)}(X-1)^{2i}\}_{0\leq i\leq k-1})\ +\  ((X-1)^{2k})\\
&=&(2^{r},2^{r-1}(X-1)^{2},2^{r-3}(X-1)^{4},...,2^{r-\nu_2((2i)!)}(X-1)^{2i},...,2^{r-\nu_2((2k-2)!)}(X-1)^{2k-2},(X-1)^{2k})\\
&&\text{and}\\
J_{r} &=& (\{2^{r-\nu_2((2i)!)}X^{2i}(X-1)^{2i}\}_{0\leq i\leq k-1}) \ +\  (X^{2k}(X-1)^{2k})\\
&=&(2^{r},2^{r-1}X^{2}(X-1)^{2},2^{r-3}X^{4}(X-1)^{4},...,\\
&&2^{r-\nu_2((2i)!)}X^{2i}(X-1)^{2i},...,2^{r-\nu_2((2k-2)!)}X^{2k-2}(X-1)^{2k-2},X^{2k}(X-1)^{2k}).\\
\end{eqnarray*}
\begin{exa}
\begin{eqnarray*}
J'_{e,1}&=&(2X,X^2),\\
J'_{e,2}&=&(4X,2X^2,X^4),\\
J'_{e,3}&=&(8X,4X^2,X^4),\\
J'_{e,4}&=&(16X,8X^2,2X^4,X^6),\\
J'_{e,5}&=&(32X,16X^2,4X^4,2X^6,X^8),\\
J'_{e,6}&=&(64X,32X^2,8X^4,4X^6,X^8).\\
\end{eqnarray*}
\end{exa}
\subsection{Generators for the full annihilator ideal}
Since $\mathbb{Z}[X]$ is noetherian, the ideals $A_e(F), A_o(F),
A(F)$ and $A_t(F)$ are finitely generated.  Under certain
conditions we give a set of generators.\\
All $f(X) \in A_t(F)$ have to vanish at $X=0 \in W_t(F)$. This
implies that the constant term of $f(X)$ is zero in the case of a
real field $F$ and that the constant term of $f(X)$ is a multiple of
$2^r$ in the case of a non real field  $F$ with level
$s(F)=2^{r-1}$. So, $A_t(F)=A_t(F) \bigcap X\mathbb{Z}[X]$ in the
real case and $A_t(F)=(A_t(F) \bigcap X\mathbb{Z}[X]) +
2^r\mathbb{Z}[X]$ in the nonreal case where the level is
$s(F)=2^{r-1}$. From now on we will study $A'_t(F) :=
A_t(F)\bigcap X\mathbb{Z}[X].$\\
 The following lemma is, in the non-real case,
proved in \cite{OVG}.  We will give a more general proof, using
Stirling numbers of the second kind.

\begin{lem}\label{incl_lemma}
Let $F$ be a field for which $2^rW_t(F)=0$ then
\ \\
\[J'_{e,r} \subseteq A'_t(F).\]
\end{lem}
\begin{proof}
\ \\
For the generator $f(X)=2^rX \in J'_{e,r}$ we have
$f(\varphi)=2^r\varphi = 0$ for all $\varphi\in W_t(F)$.\\
Let $f(X)=2^{r-\nu_2((2i)!)}X^{2i}$ for $1\leq i<k$ be one of the
other generators of $J'_{e,r}$ and $\phi \simeq \langle
a_1,a_2,\ldots , a_n\rangle$ an arbitrary even-dimensional element
of $W_t(F)$.\\ By Corollary \ref{coro}
\begin{align*}
f(\phi)&=\ \left(\sum_{q=1}^{2i}2^{q-1}1!\ S(q,1){2i\choose q}(-n)^{2i-q}2^{r-\nu_2((2i)!)}\right)(\bot_{i=1}^n\langle\!\langle a_i\rangle\!\rangle)\\
&\bot\ \left(\sum_{q=2}^{2i}2^{q-2}2!\ S(q,2){2i\choose q}(-n)^{2i-q}2^{r-\nu_2((2i)!)}\right)(\bot_{i<j}^{n}\langle\!\langle a_i,a_j \rangle\!\rangle)\\
&\bot\ \ldots\\
&\bot\ \left(\sum_{q=n}^{2i}2^{q-n}n!\ S(q,n){2i\choose
q}(-n)^{2i-q}2^{r-\nu_2((2i)!)}\right)\langle\!\langle
a_1,\ldots,a_n\rangle\!\rangle.
\end{align*}
For all $j\leq q$ we have
\begin{align*}
&\nu_2 \left(2^{q-j}j!\ S(q,j){2i\choose q}(-n)^{2i-q}2^{r-\nu_2((2i)!)}\right)\\
&\qquad\qquad=q-j+j-d(j)+\nu_2(S(q,j))+d(q)+d(2i-q)-d(2i)\\
&\qquad\qquad\qquad+(2i-q)\nu_2(n)+r-2i+d(2i)\\
&\qquad\qquad\qquad\qquad\text{(by Kummer and Legendre)}\\
&\\
&\qquad\qquad\geq q-d(j)+d(j)-d(q)+d(q)+d(2i-q)\\
&\qquad\qquad\qquad+(2i-q)\nu_2(n)+r-2i\\
&\qquad\qquad\qquad\qquad\text{(by Corollary \ref{new_lower_bound})}\\
&\\
&\qquad\qquad\geq q+d(2i-q)+(2i-q)+r-2i\\
&\qquad\qquad\qquad\qquad\text{(since $n$ is even)}\\
&\\
&\qquad\qquad\geq r
\end{align*}
Since $2^rW_t(F)=0$ it follows that
\[f(\phi)=0\] or equivalently that
\[f(X)\in A'_t(F).\]
A similar argument holds for the generator $f(X)=X^{2k}$, using the
fact that $r\leq \nu_2((2k)!)$.
\end{proof}
\ \\
This brings us to the main result of this paper:\\
\begin{thm}
Let $F$ be a field such that $2^rW_t(F)=0$ and
$2^{r-1}(I_t(F))^{2k-1}\neq 0$ with $k$ uniquely determined by
\,$\nu_2((2k-2)!)< r \leq \nu_2((2k)!)$\,. Then
\ \\
\ \\
\[J'_{e,r} = A'_t(F).\]
\end{thm}
%%%%%%%%%%%%%%%%%%%%%%%%%%%%%%%%%%%%%%%%%%%%%%%%%%%%%%%%%%%%%%%%%%%%%%%%%%%%%%%%%%%%%%%%%%%%%%%%%%%%%%%%%%%%%%%%%%
%                                                                                                                %
%  proof                                                                                                 %
%                                                                                                                %
%%%%%%%%%%%%%%%%%%%%%%%%%%%%%%%%%%%%%%%%%%%%%%%%%%%%%%%%%%%%%%%%%%%%%%%%%%%%%%%%%%%%%%%%%%%%%%%%%%%%%%%%%%%%%%%%%%
\begin{proof}
\ \\
\ \\
Let $F$ be a field such that $2^rW_t(F)=0$ and
$2^{r-1}(I_t(F))^{2k-1}\neq 0$.  Let $k$ be the unique natural
number such that $\nu_2((2k-2)!) < r \leq \nu_2((2k)!)$.\\
\\
Since $X^{2k} \in J'_{e,r}$ annihilates every even-dimensional
torsion quadratic form, we will have to prove that every polynomial
of degree $2k-1$,\\
\[f(X)=c_{2k-1}X^{2k-1}+...+c_1X,\text{ with }c_i \in \mathbb{Z},\]\\
annihilating every even-dimensional torsion quadratic form, lies in $J'_{e,r}.$\\
\\
$I_t(F)$ is generated by the elements $\langle\!\langle a
\rangle\!\rangle \in I_t(F)$ . The condition on the power of the
fundamental ideal implies the existence of elements
$a_{1},\ldots,a_{2k-1} \in F^{\ast}$ such that the form
$2^{r-1}\langle\!\langle a_1,\ldots,a_{2k-1}\rangle\!\rangle$ is
not zero.\\
Fix such elements.\\
\vskip 5pt We will evaluate the polynomial $f(X)$ in the
even-dimensional torsion quadratic forms of shape
\[ \bot_{i=1}^{n}\langle\!\langle a_{\sigma(i)}\rangle\!\rangle\]
where $\ 1\leq n\leq 2k-1,\ \sigma \in S_{2k-1}$.\\
\\
By Lemma \ref{equations} we get the following set of equations in the Witt ring:\\
\begin{equation}\label{vanishing_eq}
\left(\sum_{q=n}^{2k-1}2^{q-n}n!\ S(q,n)c_q\right)\langle\!\langle
a_1,\ldots a_n \rangle\!\rangle=0,\ 1 \leq n \leq 2k-1.
\end{equation}
\ \\
For $\ n=2k-1$, this becomes
\begin{align*}
0=& (2k-1)!\ S(2k-1,2k-1)c_{2k-1}\langle\!\langle
a_1,\ldots,a_{2k-1}\rangle\!\rangle\\
=&(2k-1)!\ c_{2k-1}\langle\!\langle
a_1,\ldots,a_{2k-1}\rangle\!\rangle.
\end{align*}
Since
\[2^{r-1}\langle\!\langle
a_1,\ldots,a_{2k-1}\rangle\!\rangle\neq 0,\] and $2^rW_t(F)=0$ it
follows that \[(2k-1)!\ c_{2k-1}=b_{2k-1}2^r\quad\text{for some}\
 b_{2k-1} \in \mathbb{Z},\] and since
$\nu_2((2k-1)!)=\nu_2((2k-2)!)<r$ that
\begin{equation*}
c_{2k-1}=2^{r - \nu_2((2k-1)!)}b_{2k-1}\quad\text{for some}\
b_{2k-1} \in \mathbb{Z}.
\end{equation*}
\\
For $\ n=2k-2\ $, (\ref{vanishing_eq}) becomes
\begin{equation}\label{second_eq}
((2k-2)!\ c_{2k-2}+2(2k-2)!\ S(2k-1,2k-2)c_{2k-1})\langle\!\langle
a_1,\ldots,a_{2k-1}\rangle\!\rangle = 0.
\end{equation}
The second term $2(2k-2)!S(2k-1,2k-2)c_{2k-1}\langle\!\langle
a_1,\ldots,a_{2k-2}\rangle\!\rangle$ vanishes since
\begin{eqnarray*}
\nu_2(2(2k-2)!\ S(2k-1,2k-2)c_{2k-1})&=&1+\nu_2((2k-2)!)\\
&&\ +\nu_2(S(2k-1,2k-2))+\nu_2(c_{2k-1})\\
&=&1+2k-2-d(2k-2)\\
&&\ +\nu_2(S(2k-1,2k-2))+\nu_2(c_{2k-1})\\
&&\qquad\qquad\qquad\qquad\qquad\qquad\text{(by Legendre)}\\
&\geq & 2k-1-d(2k-2)\\
&&\ +d(2k-2)-d(2k-1)+\nu_2(c_{2k-1})\\
&&\qquad\qquad\qquad\qquad\qquad\qquad\text{(by Corollary \ref{new_lower_bound})}\\
&\geq & 2k-1-d(2k-1)+r-(2k-1)+d(2k-1)\\
&=&r\\
\end{eqnarray*}
and $2^rW_t(F)=0.$\\
\ \\
Equation (\ref{second_eq}) is thus equivalent to\\
\begin{equation*}
(2k-2)!\ c_{2k-2}\langle\!\langle a_1,\ldots,a_{2k-2}
\rangle\!\rangle=0
\end{equation*}
and it follows, since $\nu_2((2k-2)!)=\nu_2((2k-1)!)<r$ that
\begin{equation*}
c_{2k-2}=2^{r - \nu_2((2k-2)!)}b_{2k-2}\quad\text{ for some }\
b_{2k-2} \in \mathbb{Z}.
\end{equation*}
\vspace*{2ex} \ \\
Using the same technique and observing that
\begin{align*}
\nu_2(2^{q-n}n!\ S(q,n)c_q)&=q-n+\nu_2(n!)+\nu_2(S(q,p))+\nu_2(c_q)\\
&=q-d(n)+\nu_2(S(q,n))+\nu_2(c_q)\\
&\qquad\qquad\qquad\text{(by Legendre)}\\
&\geq q-d(n)+d(n)-d(q)+r-q+d(q)\\
&\qquad\qquad\qquad\text{(by Corollary \ref{new_lower_bound})}\\
&=r\\
\end{align*}
for all $n<q$ and that $\nu_2(n!)<r$ for all $1\leq n\leq 2k-1$,
 the set of equations (\ref{vanishing_eq})
 is equivalent to the set of equations
\begin{equation*}
n!\ c_n\langle\!\langle a_1,\ldots,a_n\rangle\!\rangle=0\;\text{,
where } n=1,\ldots ,2k-1.
\end{equation*}
\ \\
The solutions are,
\[c_{n} = 2^{r-\nu_2(n!)}b_{n} \quad\text{ for some }\ b_{n}\in\mathbb{Z}.\]\\
\vspace*{3ex} We can thus rewrite,
\[f(X) = c_{2k-1}X^{2k-1}+...+c_1X\]
as
\begin{eqnarray*}
f(X)&=&2^{r-\nu_2((2k-1)!)}b_{2k-1}X^{2k-1}+2^{r-\nu_2((2k-2)!)}b_{2k-2}X^{2k-2}+\\
&&\ \ldots+2^{r-\nu_2((2j+1)!)}b_{2j+1}X^{2j+1}+b_{2j}2^{r-\nu_2((2j)!)}X^{2j}+\ldots+2^{r-\nu_2(1!)}b_1X\\
&=&2^{r-\nu_2((2k-2)!)}X^{2k-2}(b_{2k-1}X+b_{2k-2})+\ldots+2^{r-\nu_2((2j)!)}X^{2j}(b_{2j+1}X+b_{2j})+\ldots+2^rb_1X\\
\end{eqnarray*}
\vspace*{2ex}
or equivalently,\\
\[ f(X) \in J'_{e,r}\] i.e.
\ \\
\[ A'_t(F) \subset J'_{e,r}\] and using the other inclusion of the previous lemma
\ \\
\[ A'_t(F) = J'_{e,r}.\]
\end{proof}

\begin{cor}\label{coro_conj}
Let $F$ be a nonreal field of finite level $s(F)=2^{r-1}$ such
that $s(F)(I(F))^{2k-1}\neq0$, where $k=k(r)$ is uniquely
determined by $\nu_2((2k-2)!)<r\leq \nu_2((2k)!)$. Then
\[ A(F) = J_{r}.\]
\end{cor}
\begin{proof}
Since $s(F)=2^{r-1}$, we have that $2^rW(F)=0$.  Moreover, $F$ is
nonreal and hence $I_t(F) = I(F)$. The non-vanishing condition on
the power of the fundamental ideal implies that
\[ A_e(F) = A'_t(F) + (2^r) = J'_{e,r} + (2^r) = J_{e,r}.\]
If $\phi$ is an odd-dimensional form in $W(F)$, then $\phi \perp
-\,\langle 1 \rangle$ is an even-dimensional form in $W(F)$.  This
implies that
\[ A_o(F) = J_{o,r}.\]

\noindent $J_{e,r}$ and $J_{o,r}$ are comaximal ideals, since they
contain the comaximal ideals $(X^{2k})$ and $((X-1)^{2k})$
respectively.  So, we get
\[A(F) = A_e(F) \cap A_o(F) = J_{e,r} \cap J_{o,r} = J_{e,r} \cdot J_{o,r} = J_{r}.\]
\end{proof}

\begin{cor}\label{coro_real}
Let $F$ be a real field of finite height $h(F)=2^{r}$ such that
$\frac{1}{2}h(F)(W_t(F))^{2k-1}\neq0$, where $k=k(r)$ is uniquely
determined by $\nu_2((2k-2)!)<r\leq \nu_2((2k)!)$. Then
\[ A_t(F) = J_{e,r}.\]
\end{cor}
\begin{proof}
By the definition of the height $h(F)=2^{r}$, we have that
$2^rW_{t}(F)=0$. Since $F$ is a real field, we have $I_t(F)=W_t(F)$.
The non-vanishing condition on the power of the torsion ideal
implies that
\[ A_t(F) = A'_t(F) + (2^r) = J'_{e,r} + (2^r) = J_{e,r}.\]
\end{proof}
\begin{rem}
In \cite{OVG}, examples are given of non-real fields $F$ satisfying
the conditions of Corollary \ref{coro_conj}. Starting with a field
$F$ of level $2^{r-1}$, the purely transcendental extension of
transcendence degree $2k-1$, $K=F(X_1,\ldots,X_{2k-1})$ satisfies
the conditions, since the form $2^{r-1}\langle
1,X_1\rangle\otimes\ldots\otimes\langle 1,X_{2k-1}\rangle \in
s(F)(I(F))^{2k-1}$  is anisotropic
over $K$.\\
Real fields $F$ of arbitrary height $h(F)=2^r$ satisfying the conditions of Corollary \ref{coro_real} can be found
using Merkurjev's method of iterated function fields (for a description we refer the reader to \cite{DH}). We will
sketch the construction of such a real field, the details can be found in the author's thesis ( UCD - 2006). Let $k$ be
a real field and consider the purely transcendental extension of degree $(2k-1)2^r$,
$K=k(X_{1,1},\ldots,X_{1,2^r},X_{2,1},\ldots,X_{2k-1,2^r})$.  Put $\pi_i := X_{i,1}^2+\ldots+X_{i,2^r}^2$ and consider
the form $\varphi = 2^{r-1}\langle 1, -\pi_1\rangle\otimes\ldots\otimes\langle 1,-\pi_{2k-1}\rangle.$ One can show that
$\varphi$ is anisotropic over $K(2^r)$, the function field of the quadric $X_1^2+X_2^2+\ldots+X_{2^r}^2=0$ over $K$.
For a real extension $L$ of $K$ such that $\varphi$ is anisotropic over $L(2^r)$ and for any $a \in \Sigma L^2$ one
shows that $\varphi$ is anisotropic over ${L(\psi)}(2^r)$, with $\psi=2^r\langle 1,-a\rangle$. We construct a tower of
fields $K_0=K\subset K_1\subset K_2\subset \ldots$ by defining $K_{i+1}$ as the free compositum of all function fields
$K_i(\psi)$ over $K_i$ with $\psi \in \{2^r \langle 1, -b\rangle\ |\ b \in \Sigma K_i^2\}.$ Let $F=\bigcup_{i=0}^\infty
K_i$, then $F$ is a real field with Pythagoras number $p(F)=2^r$ such that $\varphi_{F} \in 2^{r-1}(W_t(F))^{2k-1}$.
Since the height $h(F)$ of a field $F$ is the smallest $2$-power greater than or equal to the pythagoras number $p(F)$
we have constructed a field satisfying the conditions of Corollary \ref{coro_real}.
\end{rem}
\begin{rem}
One can show that, for a field $F$, satisfying $2^rW(F)=0$, but not
satisfying the non-vanishing condition $2^{r-1}(I_t(F))^{2k-1}\neq
0$, the torsion annihilator ideal $A_t(F)$ always differs from
$J_{e,r}.$ To show this, one proceeds as follows. Using Corollary
\ref{coro} and the lower bound on the $2$-adic order of Stirling
numbers of the second kind, one shows that $q_l(\varphi) \in
2^{\nu_2(l!)}I^l$ for all $\varphi \in I$, with
$q_l(X)=X(X-2)(X-4)\ldots(X-2(l-1)) \in \mathbb{Z}[X].$ Choosing the
appropriate index $l=l(r)$, $q_l(X)$ annihilates all even torsion
forms and doesn't belong to $J_{e,r}.$\\
A set of generators for the ideal $A_t(F)$ in the other cases is, in
general, not known.
\end{rem}
%%%%%%%%%%%%%%%%%%%%%%%%%%%%%%%%%%%%%%%%%%%%%%%%%%%%%%%%%%%%%%%%%%%%%%%%%%%%%%%%%%%%%%%%%%%%%%%%%%%%%%%%%%%%%%%%%%
%                                                                                                                %
%  THE BIBLIOGRAPHY                                                                                              %
%                                                                                                                %
%%%%%%%%%%%%%%%%%%%%%%%%%%%%%%%%%%%%%%%%%%%%%%%%%%%%%%%%%%%%%%%%%%%%%%%%%%%%%%%%%%%%%%%%%%%%%%%%%%%%%%%%%%%%%%%%%%
%\begin{thebibliography}{999}
%\bibitem{ARA} {\bf J. Arason, R. Elman},
%Powers of the fundamental ideal in the Witt ring, {\em J. Algebra}
%{\bf239},2001, \mbox{pp. 150--160}.

\vspace*{4ex} \noindent
Stefan De Wannemacker \\
School of Mathematical Sciences \\
University College Dublin \\
Belfield, Dublin \\
Ireland \\
e-mail: {\tt sdwannem@maths.ucd.ie}
\end{document}